\documentclass{article}

\usepackage{amsmath,amsfonts,amssymb,amsthm,tikz}

\usepackage[utf8]{inputenc}
\usepackage{cite}

\newcommand{\norm}{{\mathcal{N}}}

\newcommand{\calP}{{\mathcal{P}}}

\newcommand{\C}{{\mathbb{C}}}
\newcommand{\F}{{\mathbb{F}}}

\newcommand{\Q}{{\mathbb{Q}}}

\newcommand{\height}{\mathrm{h}}

\newcommand{\gerp}{{\mathfrak{p}}}

\DeclareMathOperator{\logast}{log^\ast\!}

\newcommand{\ph}{\varphi}

    \makeatletter
    \let\@fnsymbol\@alph
    \makeatother

\title{Stewart's Theorem revisited: suppressing the norm $\pm 1$ hypothesis}
\author{Haojie Hong}

\newtheorem{theorem}{Theorem}[section]
\newtheorem{proposition}[theorem]{Proposition}

\numberwithin{equation}{section}

\makeatletter

\renewcommand*\l@section[2]{%
  \ifnum \c@tocdepth >\z@
    \addpenalty\@secpenalty
    \addvspace{0.2em \@plus\p@}%
    \setlength\@tempdima{1.5em}%
    \begingroup
      \parindent \z@ \rightskip \@pnumwidth
      \parfillskip -\@pnumwidth
      \leavevmode \bfseries
      \advance\leftskip\@tempdima
      \hskip -\leftskip
      #1\nobreak\hfil \nobreak\hb@xt@\@pnumwidth{\hss #2}\par
    \endgroup
  \fi}

\makeatother

\begin{document}

\hfuzz 4.3pt

\maketitle

\begin{abstract}
Let $\gamma$ be an algebraic number of degree $2$ and not a root of unity. In this note we show that there exists a prime ideal $\gerp$ of $\Q(\gamma)$ satisfying $\nu_\gerp(\gamma^n-1)\ge 1$, such that the rational prime $p$ underlying $\gerp$ grows quicker than $n$.
\end{abstract}

{\footnotesize

\tableofcontents

}

\section{Introduction}

Let $P(m)$ denote the largest prime factor of integer $m$, with the convention $P(0)=P(\pm 1)=1$. For any integer $n$, we denote the $n$-th cyclotomic polynomial in $x$ by $\Phi_n(x)$ as usual. 

Schinzel\cite{SC62} asked if there exists any integers $a, b$ with $ab\neq\pm 2c^2, \pm c^h(h\ge 2)$ such that $P(a^n-b^n)>2n$ for all sufficiently large $n$. Erd\H{o}s\cite{Er} conjectured that $P(2^n-1)$ grows quicker than $n$.

Let $u_n$ be the $n$th term of a Lucas sequence. In 2013, Stewart\cite{St13} gave a lower bound of the largest prime factor of $u_n$, which is of the form \\$n\exp(\log n/104\log\log n)$. What Stewart actually proved is the following, see\cite[Theorem 1.1]{BHS21}.
\begin{theorem}\label{Steor}
Let $\gamma$ be a non-zero algebraic number, not a root of unity. Denote $\omega({\gamma})$ the number of primes $\gerp$ of the field $K=\Q(\gamma)$ with the property $\nu_\gerp(\gamma)\neq 0$. Let $P$ be the biggest element of the set\[
\{\text{$p$: $p$ is a rational prime lying below a prime $\gerp$ of $K$, with $\nu_\gerp(\Phi_n(\gamma))\ge 1$}\}.
\]
If $\gamma$ satisfies one of the following conditions:
\begin{itemize}
\item
$\gamma\in\Q$,

\item $[\Q(\gamma): \Q]=2$ and $\norm\gamma=\pm 1$.
\end{itemize}
There exists a $n_0$, which is effectively computable in terms of $\omega({\gamma})$ and the discriminant of $K$, such that, for all $n>n_0$, \[
P>n\exp\left(\log n/104\log\log n\right).
\]
\end{theorem}
Using this result, Stewart answered questions of Schinzel and Erd\H{o}s in a wider context, see \cite{St13} for more historical details. 

A totally explicit expression of $n_0$ in the above theorem was given in \cite{BHS21}, which also showd that $n_0$ depends only on the field $K=\Q(\gamma)$ but not on $\omega(\gamma)$.

Let $q$ and $a$ be integers such that $q\ge 2$ and $|a|<2\sqrt{q}$. Assume $\alpha$ and $\bar{\alpha}$ are the roots of $x^2-ax+q$. In \cite{BHL21}, the authors concentrate on big prime factors of $\#E(\F_q)=q-a+1=(\alpha-1)(\bar{\alpha}-1)$, the order of group of $\F_q$-points on a certain elliptic curve $E$. A Stewart-type result was proved for recurrent sequences of order $4$ rather than Lucas sequence.

This article is motivated by \cite{BHS21} and \cite{BHL21}, we prove the following theorem.
\begin{theorem}\label{deg2}
Suppose $\gamma$ is an algebraic number of degree $2$ and not a root of unity. Set $n_0=\exp\exp(\max\{10^{10}, 3|D_K|\})$. Let $n$ be a positive integer satisfying $n\ge n_0$. There exists a prime ideal $\gerp$ of $K=\Q(\gamma)$ such that $\nu_\gerp(\gamma^n-1)\ge 1$ and the underlying rational prime $p$ of $\gerp$ satisfies\[
p\ge n\exp\left(0.0001\frac{\log n}{\log\log n}\right).
\]
\end{theorem} 

Note that this theorem suppresses the assumption $\norm\gamma=\pm 1$ when $[\Q(\gamma): \Q]=2$ in Theorem \ref{Steor}, and it can also be seen as a generalization of Schinzel's question and Erd\H{o}s' conjecture. The proof uses ingredients from \cite{BHS21} and \cite{BHL21}, all of which rely heavily on lower bound for $p$-adic logarithmic form.

\section{Preliminary results}
\subsection{Notation}
Denote by $\log^+=\max\{\log, 0\}$, $\log^-=\min\{\log ,0\}$, $\logast=\max\{\log ,1\}$.

Let $K$ be a number field of degree $d$. We denote by $D_K$ the discriminant of $K$.

Suppose $\gamma\in K$, $\height(\gamma)$ denotes the usual absolute logarithmic height of $\gamma$:\[
\height(\gamma)=[K:\Q]^{-1}\sum_{v\in M_K}d_v\log^+|\gamma|_v,
\]
where  $d_v$ denotes the local degree. The places $v\in M_K$ are normalized to extend standard places of $\Q$.

Let $\sigma: K\hookrightarrow\C$ be an arbitrary complex embedding of $K$, $\gerp$ be a prime ideal of the ring of integers $\mathcal{O}_K$. The following formulas are immediate consequences of the above definition:
\begin{equation}
\height(\gamma)=\frac{1}{d}\left(\sum_{\sigma: K\hookrightarrow\C}\log^+|\gamma^\sigma|+\sum_\gerp\max\{0, -\nu_\gerp(\gamma)\}\log\norm\gerp\right).
\end{equation}
\begin{equation}\label{hei}
\height(\gamma)=\frac{1}{d}\left(\sum_{\sigma: K\hookrightarrow\C}-\log^-|\gamma^\sigma|+\sum_\gerp\max\{0, \nu_\gerp(\gamma)\}\log\norm\gerp\right).
\end{equation}

\subsection{Uniform explicit version of Stewart's theorem}
The following two theorems go back to Stewart, see \cite[Lemma 4.3]{St13}, but in present form they are Theorem 1.4 and Theorem 1.5 of \cite{BHS21}.
\begin{theorem}
\label{thordrat}
Let~$\gamma$ be a non-zero algebraic number of degree~$d$, not a root of unity.  
Set
${p_0=\exp(80000 d (\logast d)^2)}$.  
Then for every prime~$\gerp$ of the field ${K=\Q(\gamma)}$ whose absolute norm $\norm\gerp$ satisfies ${\norm\gerp\ge p_0}$, and every positive integer~$n$   we have 
$$
\nu_\gerp(\gamma^n-1) \le \norm\gerp\exp\left(-0.002d^{-1}\frac{\log \norm\gerp}{\log\log \norm\gerp}\right)\height(\gamma)\logast n. 
$$
\end{theorem}

\begin{theorem}
\label{thordquad}
Let~$\gamma$ be a non-zero algebraic number of degree~$2$, not a root of unity.  Assume that ${\norm\gamma=\pm1}$.  Set
$
{p_0=\exp\exp(\max\{10^8,2|D_K|\})}
$,
where $D_{K}$ is the discriminant of the quadratic field ${{K}=\Q(\gamma)}$. Then for every prime~$\gerp$ of~${K}$ with underlying rational prime ${p\ge p_0}$,  and every positive integer~$n$   we have 
\begin{equation}
\label{eordquad}
\nu_\gerp(\gamma^n-1) \le p\exp\left(-0.001\frac{\log p}{\log\log p}\right)\height(\gamma)\logast n.
\end{equation}
\end{theorem}

\subsection{Cyclotomic polynomials and primitive divisors}
The following proposition, which is about eatimates of cyclotomic polynomials, goes back to Schinzel\cite{SC74}, but in the present form, item \ref{ihephin} is \cite[Theorem 3.1]{BL21} and item \ref{iarch} is proved in \cite[Proposition 8.1]{BHS21}.
\begin{proposition}
\label{prcycpol}
\begin{enumerate}
\item
\label{ihephin}
Let~$\gamma$ be an algebraic number. Then 
$$
\height(\Phi_n(\gamma)) = \ph(n)\height(\gamma) +O_1(2^{\omega(n)} \log (\pi n)),
$$
where $A=O_1(B)$ means $|A|\le B$.
\item
\label{iarch}
Let~$\gamma$ be a complex algebraic number of degree~$d$, non-zero and not a root of unity. Then 
\begin{equation}
\label{elogphi}
\log|\Phi_n(\gamma)| \ge -10^{14}d^5\height(\gamma)\cdot 2^{\omega(n)}\logast n. 
\end{equation}
\end{enumerate}
\end{proposition}

Let~$K$ be a number field of degree~$d$ and ${\gamma\in  K^\times}$ not a root of unity. We consider the sequence ${u_n=\gamma^n-1}$. Let $\gerp$ be a prime ideal of $\mathcal{O}_K$,
We call $\gerp$ \textit{primitive divisor} of $u_n$ if 
$$
\nu_\gerp(u_n)\ge 1, \qquad \nu_\gerp(u_k)= 0 \quad (k=1,\ldots n-1). 
$$
Let us recall some basic properties of primitive divisors. Items~\ref{ipd} 
of the following proposition are well-known and easy, and item~\ref{inpd} is Lemma~4 of Schinzel~\cite{SC74}; see also \cite[Lemma~4.5]{BL21}.

\begin{proposition}
\label{pprim}
\begin{enumerate} 
\item
\label{ipd}
Let~$\gerp$ be a primitive divisor of $u_n$. Then  ${\nu_\gerp(\Phi_n(\gamma) )\ge 1}$ and ${\norm\gerp\equiv1\bmod n}$; in particular, ${\norm\gerp \ge n+1}$.
%

\item
\label{inpd}
Assume that ${n\ge 2^{d+1}}$. Let~$\gerp$ be not a primitive divisor of~$u_n$. Then ${\nu_\gerp(\Phi_n(\gamma))\le \nu_\gerp(n)}$.
\end{enumerate}
\end{proposition}

\subsection{Estimates for the arimetical functions}
Denote by $\varphi(n)$, $\omega(n)$, $\tau(n)$ the Euler's totient function, the number of distnct prime factors of $n$, the number of divisors of $n$, respectively. 

We will use the following bounds for these arithmetic functions:
\begin{align}
\label{totient}
&\varphi(n)\ge 0.5\frac{n}{\log\log n}&&(n\ge 10^{20}),\\
\label{omegan}
&\omega(n)\le 1.4\frac{\log n}{\log\log n} &&(n\ge 3),\\
\label{taun}
&\tau(n)\le \exp\left(1.1\frac{\log n}{\log\log n}\right) && (n\ge 3).
\end{align}
See \cite[Theorem~15]{RS62}, \cite[Théorème~11]{Ro83}, \cite[Theorem~1]{NR83}.

\section{Proof of Theorem \ref{deg2}}

Let $P$ be the biggest element of the set\[
\text{\{$p$: $p$ is a rational prime lying below a prime $\gerp$ of $K$, with $\nu_\gerp(\Phi_n(\gamma))\ge 1$\}.}
\]
It is sufficient to show that \begin{equation}\label{mainr}
P>n\exp\left(0.0001\frac{\log n}{\log\log n}\right)
\end{equation}
One may assume $P\le n^2$, since otherwise there is nothing to prove.

By \eqref{hei},
\begin{equation}\label{byheidef}
2\height(\Phi_n(\gamma))=-\log^-|\Phi_n(\gamma)|-\log^-|\Phi_n(\gamma^\sigma)|+\sum_\gerp\max\{0,\nu_\gerp(\Phi_n(\gamma))\}\log\norm\gerp,
\end{equation}
where $\sigma$ is the generator of $\mathrm{Gal}(K/\Q)$.

We use item \ref{iarch} of Proposition \ref{prcycpol}, \begin{equation}\label{estlog}
-\log^-|\Phi_n(\gamma)| -\log^-|\Phi_n(\gamma^\sigma)|\le 2^6\cdot 10^{14}\height(\gamma)\cdot 2^{\omega(n)}\log n.
\end{equation}
We split the sum in \eqref{byheidef}:\[
\sum_\gerp\max\{0,\nu_\gerp(\Phi_n(\gamma))\}\log\norm\gerp=\sum_{\genfrac{}{}{0pt}{}{\text{$\gerp$ primi-}}{\text{tive}}}+\sum_{\genfrac{}{}{0pt}{}{\text{$\gerp$ non-}}{\text{primitive}}}
=\Sigma_{\text{p}}+\Sigma_{\text{np}},
\]
where $\gerp$ primitive, $\gerp$ non-primitive means those prime ideals $\gerp$ which are primitive, non-primitive divisors of $\gamma^n-1$, respectively.
By item \ref{inpd} of Proposition \ref{pprim},\begin{equation}\label{estnonprim}
\Sigma_{\text{np}}\le \sum_\gerp\nu_\gerp(n)\log\norm\gerp\le 2\log n.
\end{equation}
Thus 
\begin{equation}\label{upboundhe}
\height(\Phi_n(\gamma))\le10^{16}\height(\gamma)\cdot 2^{\omega(n)}\log n+\Sigma_{\text{p}}/2+\log n.
\end{equation}
On the other hand, by item \ref{ihephin} of Proposition \ref{prcycpol},
\begin{equation}\label{lowboundhe}
\height(\Phi_n(\gamma))\ge\varphi(n)\height(\gamma)-2^{\omega(n)}\log(\pi n).
\end{equation}
Combining \eqref{upboundhe} and \eqref{lowboundhe}, we have
\begin{equation}\label{sumplowerb}
\Sigma_{\text{p}}/2\ge\varphi(n)\height(\gamma)-2^{\omega(n)}\log(\pi n)-10^{16}\height(\gamma)2^{\omega(n)}\log n-\log n.
\end{equation}

Inequalities \eqref{totient}, \eqref{omegan}
and our assumption $n\ge n_0\ge 10^{10}$ 
imply that the right-hand side of \eqref{sumplowerb} is bounded from below by $0.4\varphi(n)\height(\gamma)$. 
Thus we get the lower bound of $\Sigma_{\text{p}}$
\begin{equation}\label{lowprimp}
\Sigma_{\text{p}}\ge 0.8\varphi(n)\height(\gamma).
\end{equation}
Now primes may have residue degree $1$ or $2$. Denote by
\[\Sigma_{\text{p1}}:=\sum_{\substack{\gerp~\text{primitive}\\f_\gerp=1}}\max\{0,\nu_\gerp(\Phi_n(\gamma))\}\log\norm\gerp\] and 
\[\Sigma_{\text{p2}}:=\sum_{\substack{\gerp~\text{primitive}\\f_\gerp=2}}\max\{0,\nu_\gerp(\Phi_n(\gamma))\}\log\norm\gerp.\] 
We have either
\begin{equation}\label{norm1low}
\Sigma_{\text{p1}}\ge 0.4\varphi(n)\height(\gamma),
\end{equation}
or
\begin{equation}\label{norm2low}
\Sigma_{\text{p2}}\ge 0.4\varphi(n)\height(\gamma).
\end{equation}

\subsection{Case \eqref{norm1low}}
By item \ref{ipd} of Proposition \ref{pprim}, we have $\norm\gerp=p\equiv 1\bmod n$. Since $n\ge n_0$ and $n_0=\exp\exp(\max\{10^{10}, 3|D_K|\})$ is bigger than $p_0$ in Theorem \ref{thordrat}, the underlying prime $p$ is bigger than $p_0$. 
So Theorem \ref{thordrat} applies,\[
\nu_\gerp(\Phi_n(\gamma))=\nu_\gerp(\gamma^n-1) \le p\exp\left(-0.001\frac{\log p}{\log\log p}\right)\height(\gamma)\log n. 
\]
We obtain \begin{equation}
\begin{aligned}
\Sigma_{\text{p1}}&\le \sum_{\substack{p\equiv 1\bmod n\\p\le P}}\max\{0,\nu_\gerp(\Phi_n(\gamma))\}\log p\\
&\le \pi(P; n,1)P\exp\left(-0.001\frac{\log n}{\log\log n}\right)\height(\gamma)\log n\log P,
\end{aligned}
\end{equation}
where $ \pi(P; m, a)$ counts primes $p\le x$ satisfying $p\equiv a\bmod m$. We estimate trivially $\pi(P; n,1)\le P/n$. Thus
\begin{equation}\label{norm1up}
\Sigma_{\text{p1}}\le 2P^2\exp\left(-0.001\frac{\log n}{\log\log n}\right)\height(\gamma)\frac{(\log n)^2}{n}.
\end{equation}
Combining this with \eqref{norm1low} and use \eqref{totient}, we have 
\[
P^2\ge 0.1\frac{n^2}{(\log n)^2\log\log n}\exp(0.001\frac{\log n}{\log\log n}).
\]This implies \eqref{mainr} for $n\ge n_0$.

\subsection{Case \eqref{norm2low}}
In this case, since $f_\gerp=2$, we write $p$ instead of $\gerp$.  Suppose $\sigma$ is the generator of $\mathrm{Gal}(K/\Q)$. 
For such $p$ we have $\nu_p(\gamma^n-1)=\nu_p((\gamma^\sigma)^n-1)$. Let $\beta=\gamma^\sigma/\gamma$, we  obtain the following inequalities:\[
\nu_p(\beta^n-1)\ge \nu_p((\gamma^\sigma)^n-1-(\gamma^n-1))\ge \nu_p(\gamma^n-1)\ge \nu_p(\Phi_n(\gamma)).
\]
Hence \eqref{norm2low} implies \begin{equation}\label{norm1lowb}
\sum_{p\in\calP}\nu_p(\beta^n-1)\log p\ge 0.2\varphi(n)\height(\gamma).
\end{equation}
where \[
\calP=\{p: \text{$p$ is inert in $K$ and $\nu_p(\gamma^n-1)>0$}\}.
\]
Denote $v_n=\beta^n-1$. If $\nu_p(v_n)>0$ then there exists a unique divisor $d$ of $n$ such that $p$ is primitive for $v_{n/d}$. We denote it by $d_p$. Since $\beta^n-1=\prod\limits_{d\mid n}\Phi_d(\beta)$,
\[
\nu_p(v_n)\le \nu_p(v_{n/d_p})+\sum_{\substack{m\mid n\\m\neq n/d_p}}\nu_p(\Phi_m(\beta)).
\]
By item \ref{inpd} of Proposition \ref{pprim},\[
\sum_{\substack{m\mid n\\m\neq n/d_p}}\nu_p(\Phi_m(\beta))\le \sum_{m\mid n}\nu_p(m)+\sum^7_{m=1}\nu_p(\Phi_m(\beta)).
\]
It follows that\begin{equation}\label{norm1upb}
\sum_{p\in\calP}\nu_p(\beta^n-1)\log p\le \nu_p(v_{n/d_p})\log p+\sum_{m\mid n}\log m+\sum^7_{m=1}\sum_{p\in\calP}\nu_p(\Phi_m(\beta))\log p.
\end{equation}
Trivially 
\begin{equation}\label{trivtau}
\sum_{m\mid n}\log m\le \tau(n)\log n.
\end{equation}
Notice that \[
\nu_p(\Phi_m(\beta))\le \nu_p(v_m)\le \frac{1}{2}\nu_p((\gamma^m-(\gamma^\sigma)^m)^2)
\]
and $(\gamma^m-(\gamma^\sigma)^m)^2$ is a rational integer of absolute value not exceeding \\$4(\max\{|\gamma|, |\gamma^\sigma|\})^{2m}$, we have\begin{equation}\label{vmpupb}
\sum_p\nu_p(v_m)\log p\le \log 2+m\log^+(\max\{|\gamma|, |\gamma^\sigma|\})\le \log 2+2m\height(\gamma).
\end{equation}
Hence \begin{equation}\label{7up}
\sum^7_{m=1}\sum_{p\in\calP}\nu_p(\Phi_m(\beta))\log p\le 7\log 2+56\height(\gamma).
\end{equation}
Combining \eqref{norm1lowb}\eqref{norm1upb}\eqref{trivtau}\eqref{7up}, we obtain\[
\sum_{p\in\calP}\nu_p(v_{n/d_p})\log p\ge 0.2\varphi(n)\height(\gamma)-\tau(n)\log n-56\height(\gamma)-7\log 2.
\]

\subsubsection{Big $d_p$}
Using \eqref{vmpupb},\[
\sum_{d_p\ge\tau(n)\log n}\sum_{p\in\calP}\nu_p(v_{n/d_p})\log p\le 2n\height(\gamma)\sum_{\substack{d\mid n\\d\ge\tau(n)\log n}}\frac{1}{d}+\tau(n)\log 2.
\]
The sum on the right is trivially bounded by
\[
\frac{\tau(n)}{\tau(n)\log n}=\frac{1}{\log n}.
\]
Hence\[
\sum_{d_p\ge\tau(n)\log n}\nu_p(v_{n/d_p})\log p\le \frac{2n\height(\gamma)}{\log n}+\tau(n)\log 2.
\]
Denote by $\calP'$ the subset of $\calP$ consisting of $p$ with $d_p<\tau(n)\log n$:\[
\calP'=\{p\in\calP: d_p<\tau(n)\log n\}.
\]
So we have \begin{align*}
\sum_{p\in\calP'}\nu_p(v_{n/d_p})\log p\ge &0.2\varphi(n)\height(\gamma)-\tau(n)\log n-56\height(\gamma)-7\log 2\\
&-\frac{2n\height(\gamma)}{\log n}-\tau(n)\log 2.
\end{align*}
Using \eqref{totient}\eqref{taun},  since $n\ge \exp\exp(10^{10})$, we obtain\begin{equation}\label{smalldlow}
\sum_{p\in\calP'}\nu_p(v_{n/d_p})\log p\ge 0.1\varphi(n)\height(\gamma).
\end{equation}

\subsubsection{Small $d_p: d_p<\tau(n)\log n$}
In subsections 3.2.2 and 3.2.3 of \cite{BHL21}, the authors use estimates of counting function for $S$-units to bound $\#\{d\mid n: d<\tau(n)\log n\}$ and then give a  upper bound for $\#\calP'$. One can verify that these bounds are still effective in our case, as following:
\begin{equation}\label{smalldupb}
\#\{d\mid n: d<\tau(n)\log n\}\le \exp\left(70\frac{\log n\log\log\log n}{(\log\log n)^2}\right).
\end{equation}
\begin{equation}\label{carppp}
\#\calP'\le \left(\frac{P}{n}+1\right)\exp\left(80\frac{\log n\log\log\log n}{(\log\log n)^2}\right).
\end{equation}

\subsubsection{Using Theorem \ref{thordquad}}
By item \ref{ipd} of Proposition \ref{pprim}, for $p\in\calP'$, we have $n\mid p^2-1$. Hence $p>n^{1/2}\ge n_0^{1/2}$. Since $n_0^{1/2}\ge p_0=\exp\exp(\max\{10^8,2|D_K|\})$, 
Theorem  \ref{thordquad} applies:
\begin{equation}\label{padicofb}
\begin{aligned}
\nu_p(\beta^n-1)&\le p\exp\left(-0.001\frac{\log p}{\log\log p}\right)\height(\beta)\log n\\
&\le 2P\exp\left(-0.0005\frac{\log n}{\log\log n}\right)\height(\gamma)\log n.
\end{aligned}
\end{equation}
The last inequality holds because\[
p\le P,\quad \frac{\log p}{\log\log p}\ge\frac{1}{2}\frac{\log n}{\log\log n},\quad \height(\beta)\le 2\height(\gamma).
\]
Since $\nu_p(v_{n/d_p})\le \nu_p(\beta^n-1)$, we obtain, 
\begin{equation}\label{smalldup}
\sum_{p\in\calP'}\nu_p(v_{n/d_p})\log p\le 2P\exp\left(-0.0005\frac{\log n}{\log\log n}\right)\height(\gamma)\log P\log n\#\calP'.
\end{equation}
Combining \eqref{smalldlow}\eqref{smalldup},\[
0.1\varphi(n)\height(\gamma)\le 2P\exp\left(-0.0005\frac{\log n}{\log\log n}\right)\height(\gamma)\log P\log n\#\calP'.
\]
Using \eqref{carppp} and \eqref{totient}, we obtain, for $n\ge\exp\exp(10^{10})$,
\begin{align*}
40P\left(P+n\right)\log P&\ge \frac{n^2}{\log n\log\log n}\exp\left(\frac{\log n(\log\log n-160000\log\log\log n)}{2000(\log\log n)^2}\right) \\
&\ge n^2\exp\left(0.0004\frac{\log n}{\log\log n}\right).
\end{align*}
Obviously $P\ge n$, so we have\[
80P^2\log P\ge n^2\exp\left(0.0004\frac{\log n}{\log\log n}\right),
\]
which implies \eqref{mainr} and we are done.

\paragraph{Acknowledgments}
The author thanks Yuri Bilu for checking the proof, polishing the exposition and helpful discussions. The author also acknowledges support of China Scholarship Council grant CSC202008310189.

{\footnotesize

\bibliographystyle{amsplain}
\bibliography{deg2.bib}

\providecommand{\bysame}{\leavevmode\hbox to3em{\hrulefill}\thinspace}
\providecommand{\MR}{\relax\ifhmode\unskip\space\fi MR }
\providecommand{\MRhref}[2]{%
  \href{http://www.ams.org/mathscinet-getitem?mr=#1}{#2}
}
\providecommand{\href}[2]{#2}
\begin{thebibliography}{1}

\bibitem{BHS21}
Yuri Bilu, Haojie Hong, and Sanoli Gun, \emph{Uniform explicit {S}tewart's
  theorem on prime factors of linear recurrences}, arXiv:2108.09857 (2021).
  
  
  \bibitem{BHL21}
Yuri Bilu, Haojie Hong, and Florian Luca, \emph{Big prime factors in orders of elliptic curves over finite fields}, arXiv:2112.07046 (2021).
  

  
  

\bibitem{BL21}
Yuri Bilu and Florian Luca, \emph{Binary polynomial power sums vanishing at
  roots of unity}, Acta Arith. \textbf{198} (2021), no.~2, 195--217.
  \MR{4228301}
  
  \bibitem{Er}
P. Erd\H{o}s, \emph{Some recent advances and current problems in number theory}, Lectures on
Modern Mathematics, Vol. \uppercase\expandafter{\romannumeral3}, pp. 196--244. Wiley, New York, 1965.
  

\bibitem{NR83}
J.-L. Nicolas and G.~Robin, \emph{Majorations explicites pour le nombre de
  diviseurs de {$N$}}, Canad. Math. Bull. \textbf{26} (1983), no.~4, 485--492.
  \MR{716590}

\bibitem{Ro83}
Guy Robin, \emph{Estimation de la fonction de {T}chebychef {$\theta $} sur le
  {$k$}-i\`eme nombre premier et grandes valeurs de la fonction {$\omega (n)$}
  nombre de diviseurs premiers de {$n$}}, Acta Arith. \textbf{42} (1983),
  no.~4, 367--389. \MR{736719}

\bibitem{RS62}
J.~Barkley Rosser and Lowell Schoenfeld, \emph{Approximate formulas for some
  functions of prime numbers}, Illinois J. Math. \textbf{6} (1962), 64--94.
  \MR{137689}

\bibitem{SC62}
A.~Schinzel, \emph{On primitive prime factors of $a^n-b^n$}, Proc. Cambridge Philos. Soc. \textbf{58} (1962),
  555--562. \MR{0143728}

\bibitem{SC74}
A.~Schinzel, \emph{Primitive divisors of the expression {$A^{n}-B^{n}$} in
  algebraic number fields}, J. Reine Angew. Math. \textbf{268(269)} (1974),
  27--33. \MR{344221}

\bibitem{St13}
Cameron~L. Stewart, \emph{On divisors of {L}ucas and {L}ehmer numbers}, Acta
  Math. \textbf{211} (2013), no.~2, 291--314. \MR{3143892}

\end{thebibliography}

}

\paragraph{Haojie Hong:} Institut de Mathématiques de Bordeaux, Université de Bordeaux \& CNRS, Talence, France

\end{document}